\documentclass[oneside,english]{amsart} 
\usepackage[T1]{fontenc}
\usepackage[latin9]{inputenc}
\usepackage{textcomp}
\usepackage{amsthm,amstext,amssymb}
 \usepackage[all]{xy}

\usepackage{comment}
\usepackage{pxfonts}

\makeatletter

\newcommand{\lyxmathsym}[1]{\ifmmode\begingroup\def\b@ld{bold}
  \text{\ifx\math@version\b@ld\bfseries\fi#1}\endgroup\else#1\fi}

\numberwithin{equation}{section}
\numberwithin{figure}{section}
\theoremstyle{plain}
  \newtheorem*{rem*}{Remark}
  \newtheorem{remark}{Remark}
  \newtheorem*{acknowledgement*}{Acknowledgement}
\newtheorem{theorem}[equation]{Theorem}
  \newtheorem{proposition}[equation]{Proposition}
  \newtheorem{definition}[equation]{Definition}
  \newtheorem{lemma}[equation]{Lemma}
  \newtheorem{corollary}[equation]{Corollary}
  \newtheorem{example}[equation]{Example}
  \theoremstyle{plain}
  \newtheorem*{thm*}{Theorem}
\newtheorem*{thmn*}{Theorem}
\newtheorem*{propn*}{Proposition}

\makeatother

\usepackage{babel}

\begin{document}
\title[A geometric description of the Atiyah-Hirzebruch spectral sequence]{A geometric description of the Atiyah-Hirzebruch spectral 
sequence for B-bordism}

\author{Haggai Tene}
\address{Mathematisches Institut \endgraf       
Universit\"at Heidelberg \endgraf 
Im Neuenheimer Feld 205 \endgraf 
69120 Heidelberg \endgraf 
Germany}
\email{tene@mathi.uni-heidelberg.de}
\maketitle
\begin{abstract}
In this paper we give a geometric description of the general term and the differential of the Atiyah-Hirzebruch spectral sequence for $B$-bordism. This description is given in terms of bordism classes of maps from stratifolds. We illustrate that with a computational example.  We also discuss the case of a general homology theory, where this description is given in terms of the Postnikov sections of the given theory. 

\end{abstract}

\section{Introduction}

Generalized homology theories play an important role in mathematics. Some important examples of such theories are bordism theories, including stable homotopy. Despite their relatively simple  description, computations are, in general, a hard task. One of the main tools for computations is the Atiyah-Hirzebruch Spectral Sequence (AHSS). As in other spectral sequences, the difficulty comes from computing the general terms $E^r_{p,q}$, which are certain subquotients, and the differentials $d^r_{p,q}$, which are defined using diagram chasing. The geometric nature of the theory is not completely lost in this description, but it is indirect.

In this paper we give an alternative description of the general term and the differential. Our description is geometric and makes use of stratifolds. Stratifolds, defined by Kreck (see \cite{K} and also Section \ref{sec:Stratifolds-and-Stratifold}), are certain stratified spaces generalizing smooth manifolds. In addition to the top stratum, which is a smooth manifold, stratifolds have a singular part. One can form bordism theories using stratifolds instead of manifolds. If one imposes no restriction on the stratifolds aside of compactness, and possibly a $B$-structure on their top stratum, one obtains a trivial theory. The reason is the triviality of the coefficients, since every stratifold is the boundary of its cone, which is a stratifold with boundary. To avoid this, one can impose the condition that the codimension one stratum is empty. Then the resulting theory is an ordinary homology theory. 

Another condition one might look at is that all strata of codimension $0<k<r$ are empty. For a fibration $B\to BO$, denote by $\Omega^{B(r)}$ the bordism theory of compact stratifolds with a $B$-structure on their top stratum ($B$-stratifolds), with all strata of codimension $0<k<r+2$ empty.  For $r=\infty$ we simply write $\Omega^B$ . Given a $CW$ complex $X$ (whose $k^{th}$ skeleton is denoted by $X^{k}$), for $r \geq 2$ denote
$$\hat{E}^r_{p,q}=\operatorname{Im} \left(\Omega_{p+q}^{B(q+r-2)}(X^{p})\to \Omega_{p+q}^{B(q)}(X^{p+r-1})\right).$$
This is the $r^{th}$ page of our spectral sequence. Next we describe the differential. Let $\hat{d}_{p,q}^{r}:\hat{E}^r_{p,q}\to \hat{E}^r_{p-r,q+r-1}$ be the homomorphisms induced by the map 
$$[f:\mathcal S\to X^{p}]\mapsto [g \circ f|_{\partial W }:\partial W \to X^{p-1}],$$ 
where $\mathcal S$ is a compact $B$-stratifold of dimension $p+q$ representing an element in $\Omega_{p+q}^{B(q+r-2)}(X^{p})$, $W$ is the top stratum of $S$ and $g:\partial W \to sing(\mathcal S)$ is the attaching map  which is used for gluing $W$ to $sing(\mathcal S)$, the singular part of $\mathcal S$. Note that $dim(sing(\mathcal S))\leq p-r$ so  $f$ is homotopic to a map $f'$ with  $f'\left( sing(\mathcal S)\right )\subseteq X^{p-r}$, in particular, the right side is an element in $\operatorname{Im} \left(\Omega_{p+q-1}^{B(q+2r-3)}(X^{p-r}) \to  \Omega_{p+q-1}^{B(q+r-1)}(X^{p-1}\right)$. We will show that $\hat{d}^r_{p,q}$ is well defined. 

There is a natural transformation $\varphi:\hat{E}^2_{p,q} \to H_p(X,\Omega^B_q)$ given as follows:
An element $\alpha$ in $\Omega_{p+q}^{B(q)}(X^p)$ can be represented by a map $f:\mathcal S \to X^p$  which is smooth on the preimage of the top cells, where $\mathcal S$ is a $B$-stratifold of dimension $p+q$. Then we define $\varphi(\alpha)\in H_p(X,\Omega^B_q)$ to be the homology class represented by the cycle $\Sigma_\alpha ([M_\alpha]\cdot c_\alpha)$, where each $c_\alpha$ is a $p$-cell and  $M_\alpha$ is the preimage of a regular value in $c_\alpha$ with the induced $B$-structure. We will show that it is well defined.

Our main result is the following: 

\noindent \textbf{Theorem}. 
{\it The pair $(\hat{E}^r_{p,q},\hat{d}^r_{p,q})$ is a spectral sequence and there is a natural isomorphism of spectral sequences
$$\phi:\hat{E}^r_{p,q} \to E^r_{p,q}\: ,$$ 
where the right-hand side is the standard AHSS as appears in \cite{C-F}. The natural transformation $\varphi:\hat{E}^2_{p,q} \to H_p(X,\Omega^B_q)$ is an isomorphism commuting with $\phi$ and the isomorphism $E^2_{p,q} \to H_p(X,\Omega^B_q)$ given in \cite{C-F}.}

The proof of our main theorem consists of several steps. The first is an interpretation of the AHSS in terms of the Postnikov tower, which for cohomology was done by Maunder.
Given this, the next step is a geometric interpretation of the groups occurring from the Postnikov system. The last step is the relation between the differentials and the proof that the map $\phi$ commutes with the differentials. 

In Section \ref{example} we illustrate that by giving another computation of  $\Omega_5^{fr}(\mathbb CP^{\infty})$, a result obtained by Liulevicius \cite{Lu} using the Adams Spectral Sequence. Here is a simpler example of a computation in framed bordism. Suppose $f:S^{n-1}\to S^{n-k}$ is a smooth map, which represents a non trivial element in framed bordism, where $S^{n-1}$ has the trivial framing. The mapping cone, $C_f$, has the structure of a compact framed stratifold of dimension $n$ with two strata, where the framing on the top stratum is the trivial one. For example, when $f:S^3\to S^2$ is the Hopf map then $C_f\cong \mathbb{C}P^2$. The identity map $id:C_f \to C_f$ represents an element in $\operatorname{Im}\left( \Omega^{fr(r-2)}_{n}(C_f) \to \Omega^{fr(0)}_{n}(C_f) \right)=\hat{E}^{r}_{n,0}$ 
for $r\leq k$. Its differential is represented by the map $f:S^{n-1}\to S^{n-k}$, considered as an element in $\Omega_{n-1}^{fr(r-1)}(C_f^{n-1})=\Omega_{n-1}^{fr(r-1)}(S^{n-k})$. 
When $r<k$ the differential is zero, since the mapping cylinder of $f$ is a permitted null bordism. When $r=k$ this is a non trivial element. To see this, assume that $T$ is a null bordism. Then the singular part of $T$ is of dimension at most $n-k-1$. By cellular approximation the singular part factors through the constant map, hence the map $f$ is bordant to an element in the coefficients. But this is a contradiction, since this would imply that $f$ is null bordant, since $S^{n-1}$ with its framing is null bordant.
By taking $k=2$ we get that $\hat{d}^2_{n,0}$ is, in general, non trivial for $n \geq 4$, and by taking $k=3$ we get that $\hat{d}^3_{n,0}$ is, in general, non trivial for $n \geq 5$. By analyzing all stable operations, one can see that this uniquely determines those differentials.

Steenrod's problem regarding realization of integral homology classes by maps from closed oriented manifolds was answered negatively by Thom. Replacing manifolds with stratifolds, one might ask the following:
\begin{quote}
Given a class in integral homology of a $CW$ complex, what is the minimal dimension
of the singular part in a stratifold that represents it?
\end{quote}
or how ``far" it is from being representable. As a corollary of our main theorem (Corollary \ref{cor:two filtrations agree}) it follows that this is equivalent to the question of how many steps does the element survive in the AHSS, considered as an element in $E^2_{p,0}$. 

\begin{rem*}
An alternative, more direct, proof of this fact appears in the Appendix. For this we prove a smooth approximation theorem for stratifolds, which can be useful elsewhere. 
\end{rem*}

The organization of the paper is as follows: 
In section 2 we discuss the basic properties of the Postnikov tower
of a homology theory. 
In section 3 we discuss things related to the Steenrod realization problem.
In section 4 we compare the AHSS with our spectral sequence.
In section 5 we review stratifolds, stratifold homology, and B-stratifolds. 
In section 6 we describe the Postnikov tower of a B-bordism theory using
B-stratifolds and prove our main theorem.
In section 7 we give a computational example.
In the appendix we prove a smooth approximation theorem for maps between
stratifolds. 

\section{The Postnikov tower of a homology theory }\label{section 2}
\begin{remark}
All spaces are assumed to be $CW$ complexes, and for a $CW$ complex
$X$ we denote by $X^{k}$ its $k^{th}$ skeleton. 
\end{remark}
Let $h$ be a representable generalized homology theory. One can construct the Postnikov tower, which is a sequence of homology theories $h^{(r)}$ and natural transformations between them, which fit into the following diagram

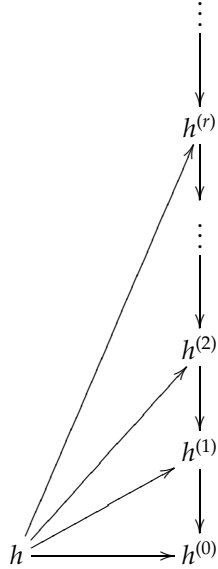
\[
\xymatrix{
& & \vdots\ar[d]\\
& & h^{(r)}\ar[d]\\
& & \vdots\ar[d]\\
& & h^{(2)}\ar[d]\\
& & h^{(1)}\ar[d]\\
h \ar[rr] \ar[rru] \ar[rruu] \ar[rruuuu]& & h^{(0)}\\
}
\]
so that it has the following properties. The theories $h^{(r)}$ have the property that the map $h_{n}\to h_{n}^{(r)}$
is an isomorphism for $n\leq r$, and $h_{n}^{(r)}$ is trivial for
$n>r$ ($h_{n}$ stands for $h_{n}(pt)$, the $n^{th}$ coefficient
group). These properties determine $h^{(r)}$ completely. For a proof
of existence and uniqueness see \cite[ Ch II, 4.13 and 4.18]{R}. 
\begin{example}
If $h$ is a connective homology theory, that is $h_{n}=0$ for $n<0$,
then $h^{(0)}$ is naturally isomorphic to homology with coefficients
in $h_{0}$ (for $CW$ complexes). If $h$ is oriented bordism, $\Omega^{SO}$,
then $h^{(0)}$ is naturally isomorphic to integral homology. Later
on we show that in the case of $\Omega^{SO}$, or other B-bordism theories, the theories $h^{(r)}$ can be given a geometric
description using B-stratifolds.\end{example}
\begin{remark}
It follows from uniqueness that for $r'\geq r$ we have $\left(h^{(r')}\right)^{(r)}=h^{(r)}$.
This implies that in all our constructions and propositions one can
replace $h$ by $ $$h^{(r')}$. 
\end{remark}
By induction, using excision, we get the following:
\begin{lemma}
\label{trivial postnikov}Let $X$ be a $CW$ complex. Then $h_{n}^{(r)}(X^{k})$
is trivial if $k+r<n$.
\end{lemma}
Again, using excision and induction, one proves the following:
\begin{lemma}
\label{lemma: relative is cellular chain}$h_{r+k}^{(r)}(X^{k},X^{k-1})\cong h_{r}\otimes C_{k}(X)$,
where $C_{k}(X)$ is the $k^{th}$ cellular chain group.
\end{lemma}
This gives a nice description of $h_{r+k}^{(r)}(X^{k})$:
\begin{lemma}
\label{lemma:cellular chains}There is a natural isomorphism  
\vspace{2mm} \\
\centerline{$h_{r+k}^{(r)}(X^{k})\cong \ker\left(h_{r}\otimes C_{k}(X)\to h_{r}\otimes C_{k-1}(X)\right).$}\end{lemma}
\begin{proof}
Look at the exact sequence of the triple $\left(X^{k},X^{k-1},X^{k-2}\right)$:

\[
0\to h_{r+k}^{(r)}(X^{k},X^{k-2})\to h_{r+k}^{(r)}(X^{k},X^{k-1})\to h_{r+k-1}^{(r)}(X^{k-1},X^{k-2}).\]
We conclude that $h_{r+k}^{(r)}(X^{k},X^{k-2})\cong \ker\left(h_{r}\otimes C_{k}(X)\to h_{r}\otimes C_{k-1}(X)\right)$.
The lemma follows from the fact that the map $h_{r+k}^{(r)}(X^{k})\to h_{r+k}^{(r)}(X^{k},X^{k-2})$
is an isomorphism using Lemma \ref{trivial postnikov}.\end{proof}
\begin{proposition}
\label{relative isomorphism for postnikov}Let $X$ be a $CW$ complex,
then the map $h_{n}(X,X^{k})\to h_{n}^{(r)}(X,X^{k})$ is an isomorphism
if $n\leq k+r+1$.\end{proposition}
\begin{proof}
We prove it for finite dimensional $CW$ complexes by induction on
the dimension. This, together with the fact that for (additive) generalized
homology theories we have \cite{M} \[
h_{*}(X)=colim(h_{*}(X^{m}))\]
will imply the statement for the case where $X$ is infinite dimensional.

The statement is trivial if $dim(X)\leq k$, since then both groups
vanish. Assume that the statement is true for every $CW$ complex $Y$ such that
$dim(Y)=m-1\geq k$, in particular for $X^{m-1}$. Using the long
exact sequence for the triple $(X^{m},X^{m-1},X^{k})$ we get the
following commutative diagram with exact rows:

 \[
\xymatrix{
h_{n+1}(X^{m},X^{m-1}) \ar[r] \ar[d]^{(1)} & h_{n}(X^{m-1},X^{k}) \ar[r] \ar[d]^{(2)}& h_{n}(X^{m},X^{k})\ar[r] \ar[d]^{(3)} & h_{n}(X^{m},X^{m-1})\ar[r] \ar[d]^{(4)} & h_{n-1}(X^{m-1},X^{k})  \ar[d]^{(5)}\\
h_{n+1}^{(r)}(X^{m},X^{m-1})\ar[r] & h_{n}^{(r)}(X^{m-1},X^{k}) \ar[r] & h_{n}^{(r)}(X^{m},X^{k}) \ar[r] & h_{n}^{(r)}(X^{m},X^{m-1}) \ar[r] & h_{n-1}^{(r)}(X^{m-1},X^{k}).
}\]
The maps $(2)$ and $(5)$ are isomorphisms by our assumption. We
have the following commutative diagram, where all horizontal maps
are isomorphisms by excision and Mayer-Vietoris:

\[
\xymatrix{
h_{j}(X^{m},X^{m-1})  \ar[r] \ar[d]^{(a)} & \oplus h_{j}(D^{m},S^{m-1}) \ar[r] \ar[d] & \oplus h_{j}(S^{m},*) \ar[r] \ar[d] & \oplus h_{j-m}(S^{0},*) \ar[d]^{(b)}\\
h_{j}^{(r)}(X^{m},X^{m-1}) \ar[r] & \oplus h_{j}^{(r)}(D^{m},S^{m-1}) \ar[r] & \oplus h_{j}^{(r)}(S^{m},*) \ar[r] & \oplus h_{j-m}^{(r)}(S^{0},*).
}\]
For $j=n$ the map $(b)$ is an isomorphism, since $j-m\leq r$, so
the same is true for $(a)$ which is $(4)$. For $j=n+1$ we either
have $\oplus h_{j-m}^{(r)}(S^{0},*)=0$ (if $n=k+r+1$ and $m=k+1$),
or else $(b)$ is an isomorphism. In any case $(b)$ is surjective
and the same is true for $(a)$ which is equal to $(1)$. Now we deduce
by the five lemma that $(3)$ is an isomorphism.
\end{proof}

We use this isomorphism to define a natural transformation \[
\Phi:h_{n}^{(r)}(X)\to h_{n-1}(X^{n-r-1})\]
as the composition \[ 
\xymatrix{
h_{n}^{(r)}(X) \ar[r] & h_{n}^{(r)}(X,X^{n-r-1}) \ar[r] & h_{n}(X,X^{n-r-1}) \ar[r] & h_{n-1}(X^{n-r-1}).
}\]
This is better seen using the following diagram:
\begin{equation}\label{diagram}
\xymatrix{
h_{n}(X) \ar[r] \ar[d] & h_{n}(X,X^{n-r-1}) \ar[r] \ar[d]^{\cong} & h_{n-1}(X^{n-r-1}) \ar[d]\\ h_{n}^{(r)}(X) \ar[r] & h_{n}^{(r)}(X,X^{n-r-1}) \ar[r] & h_{n-1}^{(r)}(X^{n-r-1}).
}
\end{equation}
The following might be a useful tool for computation:
\begin{proposition}
\label{postnikov LES}The following sequence is exact: 
\[
\xymatrix{
\cdots \ar[r]& h_n(X^{n-r-1})  \ar[r]&
h_n(X) \ar[r]&
h_n^{(r)}(X)\ar[r]^{\Phi}&
h_{n-1}(X^{n-r-1})\ar[r]&
\\
\ar[r]&
h_{n-1}^{(r)}(X^{n-r-1})\oplus h_{n-1}(X) \ar[r]&
h_{n-1}^{(r)}(X)\ar[r]&
h_{n-2}(X^{n-r-1})\ar[r]&
\cdots}
\]
%the arrow from h(r) (X) to h (Xn−r−1) is the composition of n−1 n−2
%Φ:h(r) (X) → h (Xn−r−2) with the inclusion-induced map from h (Xn−r−2) n−1 n−2 n−2
%to hn−2(Xn−r−1).
(one can prolong the sequence to the left by the sequence of the pair).
\end{proposition}
\begin{proof} 
The exactness of the sequence can be seen by diagram chasing for the Diagram \ref{diagram} (extended to the right), using
two facts: the isomorphisms \[
h_{k}(X,X^{n-r-1})\to h_{k}^{(r)}(X,X^{n-r-1})\]
 for all $k\leq n$ (Proposition \ref{relative isomorphism for postnikov})
and the fact that $h_{n}^{(r)}(X^{n-r-1})$ is trivial (Lemma \ref{trivial postnikov}).\end{proof}

\begin{corollary}
$\operatorname{Im} \left(h_{n}(X)\to h_{n}^{(r)}(X)\right)\cong \operatorname{Im} \left(h_{n}(X)\to h_{n}(X,X^{n-r-1})\right).$\end{corollary}
\begin{proof}
This follows from the fact that both kernels are equal to \[
\operatorname{Im} \left(h_{n}(X^{n-r-1})\to h_{n}(X)\right),\]
using the exact sequence in \ref{postnikov LES} and the sequence for the pair.
\end{proof}

\section{Steenrod realization problem}

Steenrod's question, which was mentioned in the introduction, can
be rephrased as asking whether the map $\Omega_{n}^{SO}(X)\to H_{n}(X)$
is surjective (here $n$ is arbitrary and the homology is with integral coefficients). Note that this map is the
localization map, using the identification $H_{n}\cong(\Omega^{SO})_{n}^{(0)}$.
Using the exact sequence in \ref{postnikov LES}, this
question is equivalent to the question whether the following map is
injective:\[
\Omega_{n-1}^{SO}(X^{n-1})\to H_{n-1}(X^{n-1})\oplus\Omega_{n-1}^{SO}(X).\]
By exactness, \[
\ker\left(\Omega_{n-1}^{SO}(X^{n-1})\to\Omega_{n-1}^{SO}(X)\right)=\operatorname{Im} \left(\Omega_{n}^{SO}(X,X^{n-1})\to\Omega_{n-1}^{SO}(X^{n-1})\right).\]
By cellular approximation, the map $\Omega_{n}^{SO}(X^{n},X^{n-1})\to\Omega_{n}^{SO}(X,X^{n-1})$
is surjective, so \[
\ker\left(\Omega_{n-1}^{SO}(X^{n-1})\to\Omega_{n-1}^{SO}(X)\right)=\operatorname{Im} \left(\Omega_{n}^{SO}(X^{n},X^{n-1})\to\Omega_{n-1}^{SO}(X^{n-1})\right).\]
$\Omega_{n}^{SO}(X^{n},X^{n-1})$ is generated by the $n$-cells,
so the right hand side is generated by the attaching maps of the $n$-cells, hence all elements are spherical.

We deduce that Steenrod's problem (in dimension $n+1$) is equivalent
to the following question.

Given a $CW$ complex $X$, does \[
\ker\left(\pi_{n}(X^{n})\to\Omega_{n}^{SO}(X^{n})\right)=\ker\left(\pi_{n}(X^{n})\to H_{n}(X^{n})\right)?\]

\begin{remark}
The equivalence is in the sense that every element which belongs to
the right side but not to the left side corresponds to a non representable
class in some $X$ with the given $n$ skeleton.
\end{remark}
Let $X$ be a $CW$ complex and $[f:M\to X]$ an element in $ $$\Omega_{*}^{SO}(X)$.
The classifying map for the stable tangent bundle of $M$ induces
a map \[
\Omega_{*}^{SO}(X)\to H_{*}(X\times BSO).\]
This map is known to be a rational isomorphism (see for example \cite{KL}
18.51).

Taking $X$ to be a point, this map is injective, i.e. one can detect
the cobordism class of $M$ by its image. It would be nice if this
was true for every space $X$. Unfortunately, this is not the case:
\begin{corollary}
The map $\Omega_{*}^{SO}(X)\to H_{*}(X\times BSO)$ need not be injective.\end{corollary}
\begin{proof}
Since we know that the answer to the Steenrod problem is negative, there is a space
$X$ and an integer $n$ such that \[
\ker\left(\pi_{n}(X^{n})\to\Omega_{n}^{SO}(X^{n})\right)\varsubsetneq \ker\left(\pi_{n}(X^{n})\to H_{n}(X^{n})\right),\]
i.e. there is a strict inclusion. Let $[f:S^{n}\to X]$ be an element
on the right side but not on the left side. Since the tangent bundle
of a sphere is stably trivial, the map $S^{n}\to X\times BSO$ factors
through $X\times*$, hence the image of $[S^{n}\to X]$ in $H_{*}(X\times BSO)$
is zero.
\end{proof}

\section{The AHSS in terms of the Postnikov tower}

For  $r \geq 2$ let $\hat{E}^r_{p,q}=\operatorname{Im} \left(h_{p+q}^{(q+r-2)}(X^{p})\to h_{p+q}^{(q)}(X^{p+r-1})\right)$, and
 $\hat{d}_{p,q}^{r}:\hat{E}^r_{p,q}\to \hat{E}^r_{p-r,q+r-1}$ be the homomorphism induced by the following diagram:
\[
\xymatrix{
 & h_{p+q}^{(q+r-2)}(X^{p}) \ar[r] \ar[d]^{\Phi} & h_{p+q}^{(q)}(X^{p+r-1})  \ar[d]^{\Phi}\\
  h_{p+q-1}^{(q+2r-3)}(X^{p-r}) \ar[r] & h_{p+q-1}^{(q+2r-3)}(X^{p-r+1}) \ar[r] & h_{p+q-1}^{(q+r-1)}(X^{p-1})  
}\]
where $\hat{E}^r_{p,q}$ is the image of the top row and $\hat{E}^r_{p-r,q+r-1}$ is the image of the composition of the bottom row. We would like to compare the pair $(\hat{E}^r_{p,q},\hat{d}^r_{p,q})$ with the AHSS. Maunder in \cite{Ma} shows that the cohomological AHSS is naturally isomorphic to the spectral sequence given by the Postnikov filtration of the cohomology theory. In 4.4 he gives a description of the general term. This can also be done in the case of homology, and one obtains the same groups as here. Maunder proves this by constructing an isomorphism between the exact couples. This way the isomorphism between the general terms is indirect, which makes is harder to identify the differentials in our case. Therefore, we give a construction of the natural isomorphism, and use it to show that the differentials agree.

\begin{theorem}\label{theorem:general term in AHSS} The pair $(\hat{E}^r_{p,q},\hat{d}^r_{p,q})$  is a spectral sequence and there is a natural isomorphism of spectral sequences
$\hat{E}^r_{p,q} \to E^r_{p,q}\: ,$
where on the right side we have the standard AHSS obtained by the exact couple. 
\end{theorem}

We start by proving the following Lemmas:
\begin{lemma}
\label{surjective maps}Let $X$ be a $CW$ complex. Then the following
maps are surjective:\\
1) $h_{n}^{(r)}(X)\to h_{n}^{(r)}(X,X^{n-r-2})$, \\
2) $h_{n}(X,X^{n-r-2})\to h_{n}^{(r)}(X,X^{n-r-2}).$ \end{lemma}
\begin{proof}
1) This follows from the long exact sequence for the pair and Lemma
\ref{trivial postnikov}.\\
2) Using the isomorphism $h_{n}(X,X^{n-r-2})\to h_{n}^{(r+1)}(X,X^{n-r-2})$
(Proposition \ref{relative isomorphism for postnikov}), it is enough
to show that the map $h_{n}^{(r+1)}(X,X^{n-r-2})\to h_{n}^{(r)}(X,X^{n-r-2})$
is surjective. We look at the following diagram:

\[
\xymatrix{
h_{n}^{(r+1)}(X,X^{n-r-2}) \ar[d]^{(1)} \ar[r]  & h_{n}^{(r+1)}(X,X^{n-r-1}) \ar[d]^{(2)} \ar[r] & h_{n-1}^{(r+1)}(X^{n-r-1},X^{n-r-2}) \ar[d]^{(3)} \\
h_{n}^{(r)}(X,X^{n-r-2}) \ar[r]^{(4)} & h_{n}^{(r)}(X,X^{n-r-1}) \ar[r]  & h_{n-1}^{(r)}(X^{n-r-1},X^{n-r-2}).}
\]
Here $(2)$ and $(3)$ are isomorphisms by Proposition \ref{relative isomorphism for postnikov}, and $(4)$ is injective since
$h_{n}^{(r)}(X^{n-r-1},X^{n-r-2})$ is trivial by Lemma \ref{trivial postnikov}. Now the lemma follows
by a diagram chase. \end{proof}
 
\begin{lemma} 
\label{butterfly lemma}Suppose we are given the following diagram:\[
\xymatrix{
A' \ar[rd]^{f'} & & \\
A \ar[r]^{f} \ar[u] & B \ar[r]^{g} & C\\
}
\]

\noindent If the bottom row is exact in $B$ then $\operatorname{Im} (f') / \operatorname{Im} (f) \cong \operatorname{Im} (g \circ f')$.
\end{lemma}

\begin{proof}
\
$ \operatorname{Im} (g \circ f') \cong \operatorname{Im} (f') / \ker{(g)} \cap \operatorname{Im} (f') = \operatorname{Im} (f') / \operatorname{Im} (f) \cap \operatorname{Im} (f') = \operatorname{Im} (f') / \operatorname{Im} (f)  $
\end{proof}

\begin{proof} (Theorem \ref{theorem:general term in AHSS})
We have the following natural isomorphism (see \cite{C-F} I,7): \[
E_{p,q}^{r}=\frac{\operatorname{Im} \left(h_{p+q}(X^{p},X^{p-r})\to h_{p+q}(X^{p},X^{p-1})\right)}{\operatorname{Im} \left(h_{p+q+1}(X^{p+r-1},X^{p})\to h_{p+q}(X^{p},X^{p-1})\right)}.\]
We use Lemma \ref{butterfly lemma} for the following diagram:
\[
\xymatrix{
h_{p+q}(X^{p},X^{p-r}) \ar[rd] \ar[rrd]^{f_1}& & \\
h_{p+q+1}(X^{p+r-1},X^{p}) \ar[r] \ar[u] & h_{p+q}(X^{p},X^{p-1}) \ar[r] & h_{p+q}(X^{p+r-1},X^{p-1}),\\
}
\]
where $f_1$ is the composition,
%We use lemma \ref{butterfly lemma} with $A'=h_{p+q}(X^{p},X^{p-r}), A=h_{p+q+1}(X^{p+r-1},X^{p}), B=h_{p+q}(X^{p},X^{p-1}), C=h_{p+q}(X^{p+r-1},X^{p-1})$. 
to conclude that:
$$E^r_{p,q} \cong \operatorname{Im} \left( f_1\right).$$
In order to simplify this expression, look at the following diagram:
\[
\xymatrix{h_{p+q}(X^{p},X^{p-r}) \ar[r]\ar[d]  \ar@/^1pc/[rr]^{f_2}
\ar[rd]^{f_1}& h_{p+q}(X^{p+r-1},X^{p-r}) \ar[r]\ar[d] & h_{p+q}^{(q)}(X^{p+r-1},X^{p-r}) \ar[d]^{(1)}\\
h_{p+q}(X^{p},X^{p-1})  \ar[r] & h_{p+q}(X^{p+r-1},X^{p-1})  \ar[r]^{(2)} & h_{p+q}^{(q)}(X^{p+r-1},X^{p-1}).
}
\]
where $f_2$ is the composition.\\
$(1)$ is injective by the long exact sequence for the triple $\left(X^{p+r-1},X^{p-1},X^{p-r}\right)$
and Lemma \ref{trivial postnikov};\\
$(2)$ is an isomorphism by Proposition \ref{relative isomorphism for postnikov}.
%Denote $f:h_{p+q}(X^{p},X^{p-r})\to h_{p+q}^{(q)}(X^{p+r-1},X^{p-r})$.

\noindent{From $(1)$ and $(2)$ it follows that $E_{p,q}^{r}\cong \operatorname{Im} \left(f_1\right) \cong \operatorname{Im} \left(f_2\right)$.
Look at the following diagram:}
\[
\xymatrix{& h_{p+q}^{(q+r-2)}(X^{p}) \ar[d]^{(2)} \ar[r] & h_{p+q}^{(q)}(X^{p+r-1})\ar[d]^{(3)} \\
h_{p+q}(X^{p},X^{p-r}) \ar[r]^(.45){(1)}  \ar@/_1pc/[rr]_{f_2}   \ar[r] & h_{p+q}^{(q+r-2)}(X^{p},X^{p-r}) \ar[r] & h_{p+q}^{(q)}(X^{p+r-1},X^{p-r})}
\]
$(1)$ and $(2)$ are surjective by Lemma \ref{surjective maps};\\
$(3)$ is an isomorphism by the long exact sequence for the pair and Lemma \ref{trivial postnikov}.\\
This implies that $\operatorname{Im}(f_2)\cong \hat{E}^r_{p,q}.$

The fact that this isomorphism commutes with the differential follows from diagram chasing.
\end{proof}

\subsection*{The $\hat{E}^{2}$ page}

The case $r=2$ has the following form:\[
\hat{E}_{p,q}^{2}=\operatorname{Im} \left(h_{p+q}^{(q)}(X^{p})\to h_{p+q}^{(q)}(X^{p+1})\right).\]
Note that by the long exact sequence for the pair, we have\[
\operatorname{Im} \left(h_{p+q}^{(q)}(X^{p})\to h_{p+q}^{(q)}(X^{p+1})\right)\cong h_{p+q}^{(q)}(X^{p})/\operatorname{Im} \left(h_{p+q+1}^{(q)}(X^{p+1},X^{p})\to h_{p+q}^{(q)}(X^{p})\right),\]
and by Lemmas \ref{lemma: relative is cellular chain} and \ref{lemma:cellular chains}
this is naturally isomorphic to\[
\ker\left(h_{q}\otimes C_{p}(X)\to h_{q}\otimes C_{p-1}(X)\right)/\operatorname{Im} \left(h_{q}\otimes C_{p+1}(X)\to h_{q}\otimes C_{p}(X)\right),\]
which is by definition $H_{p}(X,h_{q})$, as we know from the standard
presentation of the AHSS. It is not hard to see this is compatible with the isomorphism $E^2_{p,q} \to H_p(X,\Omega^B_q)$ which appears in \cite{C-F}.

\section{\label{sec:Stratifolds-and-Stratifold}Stratifolds, Stratifold
homology, and generalization to B-stratifolds theories}

Stratifolds are generalization of manifolds. They were introduced
by Kreck \cite{K} and used in order to define a bordism theory, denoted
by $SH_{*}$, which is naturally isomorphic to singular homology for
$CW$ complexes. Kreck also defined a cohomology theory using stratifolds
which is defined on the category of smooth oriented manifolds (without
boundary but not necessarily compact). It is denoted by $SH^{*}$
and is naturally isomorphic to singular cohomology.

\subsection*{Stratifolds and B-stratifolds}

Kreck defined stratifolds as spaces with a certain sheaf of functions,
called the smooth functions, fulfilling certain properties, but for
our purpose the following definition is enough (these stratifolds
are also called p-stratifolds).

Stratifolds are constructed inductively in a similar way to the way
we construct $CW$ complexes. We start with a discrete set of points
denoted by $X^{0}$ and define inductively the set of smooth functions,
which in the case of $X^{0}$ are all real functions. 

Suppose $X^{n-1}$ together with a smooth set of functions is given.
Let $W$ be an $n$-dimensional smooth manifold, {}``the $n$ stratum''
with boundary and a collar $c$, and $f$ a continuous map from the
boundary of $W$ to $X^{n-1}$. We require that $f$ will be proper
and smooth, which means that its composition with every smooth map
from $X^{n-1}$ is smooth. Define $X^{n}=X^{n-1}\cup_{f}W$. The smooth
maps on $X^{n}$ are defined to be those maps $g:X^{n}\to\mathbb{R}$
which are smooth when restricted to $X^{n-1}$ and to $W$ and such
that for some $0<\delta$ we have $gc(x,t)=gf(x)$ for all $x\in\partial W$
and $t<\delta$.

Among the examples of stratifolds are smooth manifolds, real algebraic varieties
\cite{Gri}, and the one point compactification of a smooth manifold (which is the interior of a manifold with boundary).
The cone over a stratifold and the product of two stratifolds are
again stratifolds. 

We can also define stratifolds with boundary, which are analogous
to manifolds with boundary. A main difference is that every stratifold
is the boundary of its cone, which is a stratifold with boundary. 

Given two stratifolds with boundary $(T',S')$ and $(T'',S'')$ and
an isomorphism $f:S'\rightarrow S''$, there is a well defined stratifold
structure on the space $T'\cup_{f}T''$, that is called the gluing.
On the other hand, given a smooth map $g:T\rightarrow\mathbb{R}$
such that there is a neighborhood of $0$ which consists only of regular
values then the preimages $g^{-1}((-\infty,0])=T'$ and $g^{-1}([0,\infty))=T''$
are stratifolds with boundary, and $T$ is isomorphic to the gluing
$T'\cup_{Id}T''$.

To obtain singular homology we specialize our stratifolds in the following
way: we use compact stratifolds, require that their top stratum will be
oriented and the codimension one stratum will be empty. 

One can generalize this in several ways. One way will be to require some B-structure on the top stratum, like a spin or a string structure. We will call such stratifolds B-stratifolds. Another way will be to restrict the allowed strata. In the next sections we will discuss some of these generalizations.

\begin{remark}
Regularity, a condition that is often required, is not needed
here as was noted by Kreck in his preprint {}``Pseudo Homology, p-Stratifold
Homology and Ordinary Homology''.
\end{remark}

\subsection*{Stratifold homology}

Stratifold homology was defined by Kreck in \cite{K}. We will describe
here a variant of this theory called parametrized stratifold homology,
which is naturally isomorphic to it for $CW$ complexes. In this paper
we will refer to parametrized stratifold homology just as stratifold
homology and use the same notation for it. 

\begin{definition}
Let $X$ be a topological space and $n\geq0$, define $SH_{n}(X)$
to be $\left\{ g:S\rightarrow X\right\} /\sim$, i.e., bordism classes
of maps $g:S\rightarrow X$, where $S$ is a compact oriented stratifold
of dimension $n$ and $g$ is a continuous map. We often denote the
class $[g:S\rightarrow X]$ by $[S,g]$ or by $[S\rightarrow X]$.
$SH_{n}(X)$ has a natural structure of an Abelian group, where addition
is given by disjoint union of maps and the inverse is given by reversing
the orientation. If $f:X\rightarrow Y$ is a continuous map, then
the induced map $f_{*}:SH_{n}(X)\rightarrow SH_{n}(Y)$ is given by
composition. 
\end{definition}
One constructs a boundary operator and proves (Kreck, \cite{K}, chapter 5):
\begin{theorem}
(Mayer-Vietoris) The following sequence is exact:\[
 \cdots \rightarrow SH_{n}(U\cap V)\rightarrow SH_{n}(U)\oplus SH_{n}(V)\rightarrow SH_{n}(U\cup V)\xrightarrow{\partial}SH_{n-1}(U\cap V)\rightarrow \cdots \]
where the first map is induced by inclusions and the second is the
difference of the maps induced by inclusions. 
\end{theorem}
$SH_{*}$ with the boundary operator is a homology theory. Its main property is the following (Kreck, \cite{K}, 20.1):
\begin{theorem}
There is a natural isomorphism of homology theories $\varphi:SH_{*}\to H_{*}$.
\end{theorem}
$\varphi$ is given by $\varphi_{n}([S,f])=f_{*}([S])$,
where $[S]\in H_{n}(S,\mathbb{Z})$ is the fundamental class of $S$.
\begin{remark}
One can replace stratifolds with B-stratifolds and all the construction will still work. One only needs to know how to give a B-structure to a boundary and to see that gluing along a boundary gives a B-structure. We denote these theories by $SH_{*}^{B}$ and call them B-stratifold homology. 
\end{remark}

\section{The Postnikov tower of a B-bordism theory}
 
Define the following sequence of variants of B-stratifold homology:\[
\Omega_{p}^{B (k)}(X)=[f:S\to X]/\sim,\]
 where $S$ is a compact B-stratifold which has
empty strata in codimension $<k+2$ or, equivalently, its singular
part is of dimension at most $p-k-2$. The same condition for the codimension must hold
for the bordism relation. Note that $\Omega_{p}^{B (0)}(X)=SH_{p}^{B} (X)$. In case $B=SO$ we omit the $B$ from the notation. 
Clearly, there are the following natural transformations: \[
\Omega_{p}^{B}(X)\to \cdots \to \Omega_{p}^{B (2)}(X)\to \Omega_{p}^{B (1)}(X)\to \Omega_{p}^{B (0)}(X) = SH_{p}^{B}(X).\]

The map $\Omega_{p}^{B}\to \Omega_{p}^{B (r)}$ is an isomorphism for $p\leq r$ since in this range the singular part of the cycles and bordisms in $\Omega_{n}^{B (r)}(pt)$ must be empty. $\Omega_{p}^{B (r)}$ is trivial for $p>r$ since then the cone of a cycle is an allowed null bordism. This proves the following:
\begin{proposition}\label{postnikov via stratifolds}
$\left(\Omega^{B}\right)^{(r)}\cong \Omega^{B (r)}$.
In particular, $SH^{B} \cong H( -,\Omega^{B}_{0})$.
\end{proposition}

We use this geometric description of $\left(\Omega^{B}\right)^{(r)}$
in order to describe the general term in the AHSS for $\Omega^{B}$:
Recall our definition $$\hat{E}^r_{p,q}=\operatorname{Im} \left(\Omega_{p+q}^{B(q+r-2)}(X^{p})\to \Omega_{p+q}^{B(q)}(X^{p+r-1})\right),$$
and the differentials $\hat{d}^r_{p,q}:\hat{E}^r_{p,q}\to \hat{E}^r_{p-r,q+r-1}$ given by: 
$$[f:\mathcal S\to X^{p}]\mapsto [g \circ f|_{\partial W }:\partial W \to X^{p-1}],$$ 
as noted in the introduction. We prove our main theorem:

\begin{proof} (Main Theorem)
The identification of the spectral sequences follows from Theorem \ref{theorem:general term in AHSS} and Proposition \ref{postnikov via stratifolds}. The differential $d^r_{p,q}$ is induced by the natural transformation $$\Phi:\Omega_{n}^{B (r)}(X)\to \Omega^B_{n-1}(X^{n-r-1}),$$ which was discussed in section 2. A direct translation of this natural transformation gives our differential $\hat{d}^r_{p,q}$ (This is easily seen by the fact that the inverse of the isomorphism \[
\Omega^{B}_{n}(X,X^{n-r-1}) \to \Omega^{B (r)}_{n}(X,X^{n-r-1})
\]
is given by restriction to the top stratum). Hence the isomorphism commutes with the differentials. This implies that the pair $(\hat{E}^r_{p,q},\hat{d}^r_{p,q})$ is a spectral sequence. The identification of the second page follows from the identification in the general case appearing in section 2.
\end{proof}

We can compare this description of the differential to other known descriptions. One example is the $d^2$ differential in spin bordism in the rows $q=0,1$ (Lemma 2.3.2 in \cite{Te}):

\vspace{2mm}
\centerline{$d^2_{p,0}:H_p(X,\mathbb{Z}) \to H_{p-2}(X,\mathbb{Z}/2)$}
\vspace{2mm}

\noindent is given by reduction mod $2$ composed with the dual of  $sq^2$.
\vspace{2mm}

\centerline{$d^2_{p,1}:H_p(X,\mathbb{Z}/2) \to H_{p-2}(X,\mathbb{Z}/2)$}
\vspace{2mm}

\noindent is given by the dual of  $sq^2$.

Let us describe $\hat{d}^2_{p,0}$. An element $\alpha \in H_p(X,\mathbb{Z})$ can be represented (uniquely, up to bordism) by a map from a compact $p$-dimensional spin stratifold $[S\to X]$. The differential is given by the restriction to the boundary of the top stratum, say $M^{p-1}\to X$. Since this map factors through the singular part of $S$, we can assume that the image is contained in $X^{p-2}$. Composition with the collapse map $X^{p-2}\to X^{p-2}/X^{p-3}\cong \bigvee S^{p-2}$ gives a $1$-cycle since a map $M^{p-1}\to S^{p-2}$ gives an element in $\Omega^{Spin}_1\cong \mathbb{Z}/2$. This gives a nice geometric description to the dual of $sq^2$, at least for classes which are in the image of the reduction map from integral homology.

\subsection*{The filtration in homology}
When $q=0$, we have that\[
\hat{E}_{p,0}^{r}=\operatorname{Im} \left(\Omega_{p}^{B (r-2)}(X^{p})\to SH^{B}_{p}(X^{p+r-1})\right).\]
When $r\geq2$, this is equal to\[
\hat{E}_{p,0}^{r}=\operatorname{Im} \left(\Omega_{p}^{B (r-2)}(X)\to SH^{B}_{p}(X)\right).\]
Assume that $\Omega^{B}_{0} \cong \mathbb{Z}$ so $SH^{B}_{p}(X) \cong H_{p}(X;\mathbb{Z})$, which is often the case, then we conclude:
\begin{corollary}
\label{cor:two filtrations agree}For a $CW$ complex $X$, the filtration
in singular homology given by the AHSS for B-bordism\[
E_{p,0}^{\infty}\subseteq \cdots \subseteq E_{p,0}^{4}\subseteq E_{p,0}^{3}\subseteq E_{p,0}^{2}\cong H_{p}(X,\mathbb{Z}),\]
agrees with the filtration given by all classes in homology that are
represented by maps from B-stratifolds with singular part of dimension
at most $p-r-2$.\end{corollary}

In the appendix (Corollary \ref{corollary:The-filtrations-given}) we give a more straightforward proof of this fact.

\section{An example}\label{example}

To illustrate our description, we compute the group $\Omega^{fr}_5(\mathbb C P^{\infty})$. This group was first computed by Liulevicius in \cite{Lu}, using the Adams Spectral Sequence. We start by computing  $\hat{d}^2_{4,1}:\hat{E}^2_{4,1}\to \hat{E}^2_{2,2}.$ Using Lemma 
\ref{lemma:cellular chains} we have 
$$\hat{E}^2_{4,1}=\Omega^{fr(1)}_5(\mathbb C P^2) \cong h_1\otimes C_4(\mathbb C P^{\infty}) \cong \mathbb Z /2,$$ 
$$\hat{E}^2_{2,2}=\Omega^{fr(2)}_4(\mathbb C P^1) \cong h_2\otimes C_2(\mathbb C P^{\infty}) \cong \mathbb Z /2.$$ 
We describe the generator of  $\Omega^{fr(1)}_5(\mathbb C P^2)$. Let $\mathcal S = D^4\times S^1 \cup _g \mathbb C P^1$, where $D^4$ has the standard framing, $S^1$ has the Lie group framing, and the gluing map  $g:S^3\times S^1 \to \mathbb C P^1$ is given by first projecting on $S^3$ and then the Hopf map. $\mathcal S$ is a framed stratifold of dimension $5$ with a singular part of dimension $2$. The map $f:\mathcal S \to \mathbb C P^2$ is given as follows: $f:D^4\times S^1\to D^4$ is the projection and $f:\mathbb C P^1 \to \mathbb C P^1$ is the identity.  $(\mathcal S,f)$ is an allowed cycle in $\Omega^{fr(1)}_5(\mathbb C P^2)$ since the singularity is of codimension $3$. The class $\hat{d}^2_{4,1}([\mathcal S,f])$ is given by $[S^3\times S^1\to \mathbb C P^1]$. This class is non trivial since it is a fiber bundle where all fibers are non null bordant, and hence $\hat{d}^2_{4,1}$ is an isomorphism. 

Next, we compute the differentials coming out of $\hat{E}^r_{6,0}=\Omega^{fr(0)}_6(\mathbb C P^3)\cong \mathbb Z$, by constructing a generator for this group. Consider the fiber bundle 
\[
\xymatrix{
\mathbb CP^1 \ar[r]& \mathbb CP^3 \ar[r]^{p}  &S^4.
}
\]
Set $\mathcal S=\mathbb CP^3$ as a framed stratifold the following way:  Choose a point $x \in S^4$ and set $\mathbb CP^1\cong p^{-1}(x)$ to be the singular part of $\mathcal S$. Set $D=p^{-1}(S^4 \setminus \{x\})$ to be the top stratum of $\mathcal S$ together with the framing obtained from $S^4 \setminus \{x\}$ and the trivial framing of the fiber $\mathbb CP^1$. Note that by choosing a small disk centered at x, the preimage of its complement is a framed manifold with boundary $\partial$ and we can consider $\mathcal S$ as the gluing of this manifold to $\mathbb CP^1$ along the projection map $g:\partial \to \mathbb CP^1$. The identity map $f:\mathcal S \to \mathbb CP^3$ is a generator of $\Omega^{fr(0)}_6(\mathbb C P^3)$. The fact that the singular part of $\mathcal S$ is of dimension 2 implies that this generator is mapped to zero by $\hat{d}^2_{6,0}$ and $\hat{d}^3_{6,0}$ (see Corollary 
\ref{cor:two filtrations agree}). We are left with computing $\hat{d}^4_{6,0}$. This differential is induced by the following diagram

\[
\xymatrix{
\Omega^{fr(2)}_6(\mathbb CP^3)            \ar[r] \ar[d]&            \Omega^{fr(0)}_6(\mathbb CP^4)  \cong  \mathbb Z \ar[d]\\
\mathbb Z/24  \cong     \Omega^{fr(5)}_5(\mathbb CP^1)            \ar[r] &             \Omega^{fr(3)}_5(\mathbb CP^2).
}
\]
Note that the only differential entering the $(2,3)$ position before the $4^{th}$ page is coming from $$\hat{E}^2_{4,2}=\Omega^{fr(2)}_6(\mathbb CP^2)\cong \mathbb Z/2.$$ Hence, the kernel of the bottom map in the diagram is either trivial or $12\mathbb Z/24$. We will show that if we lift a generator of $ \Omega^{fr(0)}_6(\mathbb CP^4) $ to $\Omega^{fr(2)}_6(\mathbb CP^3)$ and then map it to $\Omega^{fr(5)}_5(\mathbb CP^1)$ we get twice the generator (Note that $\Omega^{fr(5)}_5(\mathbb CP^1)\cong\Omega^{fr}_5(\mathbb CP^1)$). This will imply that $\hat{E}^5_{2,3}\cong \mathbb Z /2$. Taking $[\mathcal S \to \mathbb CP^3]$ as a generator, the element in $\Omega^{fr}_5(\mathbb CP^1)$ is given by the gluing map $[g:\partial \to \mathbb CP^1]$. Note that $g$ is a (trivial) bundle whose fiber is $S^3$, so it is enough to analyze the framing on the fiber induced by the framing on $\partial$, which itself is induced by the framing on the preimage of the complement of a disk around $x$. It is well known that the bundle $\mathbb CP^1 \to \mathbb CP^3 \xrightarrow{p} S^4$ is linear and it is obtained by the clutching function given by the generator in $\pi_3(SO(3))$, which is mapped to twice the generator in $\pi_3(SO)$. Thus the induced framing on the fiber is twice the generator in  $\Omega^{fr}_3$, and we are done.

\appendix
\section{Smooth approximation of maps between stratifolds}

We give here a more direct proof of Corollary \ref{cor:two filtrations agree}. The main tool is the approximation theorem.
\begin{lemma}
\label{lemma:Every-continuous-map}Every continuous map $f:M\to S$,
where $M$ is a manifold with boundary $\partial M$, $S$ is a stratifold,
and $f$ is smooth when restricted to $\partial M$, is homotopic
to a smooth map rel. boundary.\end{lemma}
\begin{proof}
We prove it by induction on the dimension of $S$. This is clear if
$S$ is $0-$dimensional, since then $S$ is discrete, so $f$ is
locally constant. Assume that we know that for stratifolds of dimension
$<n$, and let $f:M\to S$ be a continuous map, where $S$ is of dimension
$n.$ 

By using the collar of $\partial M$, we can assume that $f$ is smooth
on a neighbourhood of $\partial M$. 

By construction, $S$ is obtained from $\Sigma$, the singular part,
and a manifold $N$ of dimension $n$ with a boundary $\partial N$
and an (open) collar. Denote by $U$ the singular part together with
the collar, and by $V$ the interior of $N$, that is $N\setminus\partial N$.
This is an open cover of $S$, denote $f^{-1}(U)=U'$ and $f^{-1}(V)=V'$,
then this is an open cover of $M$. We can choose a smooth function
$g:M\to\mathbb{R}$ such that $g|_{U'\setminus V'}=0$ and $g|_{V'\setminus U'}=1$,
and a regular value both of $g$ and of $g|_{\partial M}$, say $0.5$.
Denote by $P$ its preimage. Then $P$ is a manifold with a boundary,
denoted by $\partial P$. By a standard approximation argument one
can show that $f$ is homotopic rel. boundary to a map $\tilde{f}:M\to S$
with the following properties:\\
1) $\tilde{f}$ is smooth in a neighbourhood of $\partial M\cup g^{-1}([0.5,1])$;\\
2) $g^{-1}([0.5,1])$ is mapped into $ $$V$;\\
3) $g^{-1}([0,0.5])$ is mapped into $ $$U$.

We are left with smoothing $g^{-1}([0,0.5))$. We can find a manifold
$M'$ of dimension $n$ with boundary $\partial M'$, embedded in
$g^{-1}([0,0.5))$ such that $\tilde{f}$ is smooth outside of $M'$
and in a neighbourhood of $\partial M'$. $M'$ is mapped to $U$,
which is smoothly homotopy equivalent to $\Sigma$. This implies that
it is homotopic to a new map, which is smooth outside of $M'$ and
in a neighbourhood of its boundary, and that the image of $M'$ is
contained in $\Sigma$. Now we use the inductive step to smooth this
map using the fact that $\Sigma$ is of dimension $<n$.\end{proof}
\begin{proposition}
Every continuous map $f:S\to S'$, where $S$ and $S'$ are stratifolds
is homotopic to a smooth map. \end{proposition}
\begin{proof}
This is proved by induction on the dimension of $S$ using Lemma \ref{lemma:Every-continuous-map}.\end{proof}
\begin{proposition}
\label{finite CW=stratifold}Any locally finite, finite dimensional countable
$CW$ complex is homotopy equivalent to a stratifold with the same set of cells. \end{proposition}
\begin{proof}
This is proved by induction on the dimension using the fact that the
attaching maps can be made smooth using Lemma \ref{lemma:Every-continuous-map}.\end{proof}
\begin{lemma}
\label{lemma:Let--be7}Let $X^{d}$ be a $CW$ complex of dimension
$d$ and $M^{m}$ a closed oriented manifold of dimension $m>d$ together
with a map $f:M^{m}\to X^{d}$. The element $(M^{m},f)$ represents
the zero element in $SH_{m}(X^{d})$, and there exists a null bordism,
which has singular part of dimension $\leq d$.\end{lemma}
\begin{proof}
The first assertion is clear since $SH_{m}(X^{d})=0$. First assume
that $X^{d}$ is a compact stratifold and $f$ is smooth. In this
case the null bordism can be taken to be the mapping cylinder $Cyl(f)$
whose boundary is $M^{m}$ and whose singular part is $X^{d}$. The
general case follows from Proposition \ref{finite CW=stratifold}
since the image of $f$ is contained in some finite subcomplex that
is homotopy equivalent to a stratifold.\end{proof}
\begin{corollary}
\label{corollary:The-filtrations-given}The filtrations given by the images
\[
\Omega_{p}(X)\to \cdots \to \Omega_{p}^{(2)}(X)\to \Omega_{p}^{(1)}(X)\to \Omega_{p}^{(0)}(X)\cong H_{p}(X,\mathbb{Z})\]
and \[
\Omega_{p}(X)\to \cdots \to\Omega_{p}(X,X^{p-3})\to\Omega_{p}(X,X^{p-2})\to H_{p}(X,X^{p-2})\cong H_{p}(X,\mathbb{Z})\]
are equal.\end{corollary}
\begin{proof}
Given $[S,f]\in \Omega_{p}^{(k)}(X)$ we can assume by cellular approximation
that $S^{p-k-2}$, the $p-k-2$ skeleton of $S$, is mapped to the
$X^{p-k-2}$, the $p-k-2$ skeleton of $X$. Denote by $(N^{p},\partial N^{p})$
the manifold we used to get the top stratum of $S$. Then there is
an induced map $(N^{p},\partial N^{p})\xrightarrow{f}(X,X^{p-k-2})$.
By the definition of the bordism relation, $[(S,\emptyset),f]=[(N^{p},\partial N^{p}),f]\in \Omega_{p}(X,X^{p-k-2})$. 

Let $[(N^{p},\partial N^{p}),f]\in\Omega_{p}(X,X^{p-k-2})$ be any
element. By Lemma \ref{lemma:Let--be7}, $(\partial N^{p},f|_{\partial N^{p}})$
represents the zero element in $SH_{m}(X^{p-k-2})$, and there exists
a null bordism, $(S,\partial N^{p})$, which has singular part of
dimension $\leq p-k-2$. Gluing $(N^{p},\partial N^{p})$ and $(S,\partial N^{p})$
we get a stratifold $N^{p}\cup_{\partial N^{p}}S$ of dimension $p$
with a singular part of dimension $\leq p-k-2$ mapped to $X$, hence
an elements in $\Omega_{p}^{(k)}(X)$, such that its image in $SH_{p}(X,X^{p-k-2})$
equals to the image of $[(N^{p},\partial N^{p}),f]$.
\end{proof}
Lemma \ref{lemma:Let--be7} and its proof have another nice
\begin{corollary}
Let $X$ be a $CW$ complex having cells only in even dimensions. Every homology class can be represented by
a map from a stratifold having only even dimensional strata. This
can be restated in the following way: let $SH_{*}^{even}$ be the
bordism theory of stratifolds with strata only in even codimensions.
Then the map $ $$SH_{*}^{even}(X)\to SH_{*}(X)$ is surjective.\end{corollary}
\begin{proof}
In this case the odd dimensional homology groups are trivial, so we
may only look at even dimensional ones. Given $\alpha\in SH_{2k}(X)$,
we represent it using a stratifold $[S,f]$, and like we did before,
we can replace the singular part of $S$ with the $2k-2$ skeleton
of $X$ without changing the homology class.\end{proof}
\begin{acknowledgement*}
I would like to thank Matthias Kreck and Diarmuid Crowley for many helpful conversations.
A large part of this work was done while I was member of Graduiertenkolleg
1150 Homotopy and Cohomology at the University of Bochum. I was supported
by Priority Research Centers Program through the National Research
Foundation of Korea (NRF) funded by the Ministry of Education, Science
and Technology (project No. 2012047640) while I was at PMI/POSTECH.

\end{acknowledgement*}

\end{document}